\documentclass[10pt]{article}
\usepackage{latexsym}
\usepackage{amssymb}

\newtheorem{teo}{THEOREM}[section]
\newtheorem{prop}[teo]{PROPOSITION}
\newtheorem{lem}[teo] {LEMMA}
\newtheorem{ejem}[teo]{Example}
\newtheorem{obser}[teo]{Remark}
\newtheorem{defi}[teo]{DEFINITION}
\newtheorem{coro}[teo]{COROLLARY}
\newenvironment{Demo}{\noindent \sf Proof. \rm}

\def\Hom{\mathop{\rm Hom}\nolimits}

\def\qed{\hfill \mbox{$\square$}}
\def\End{\mathop{\rm End}\nolimits}

\def\Ext{\mathop{\rm Ext}\nolimits}
\def\Tor{\mathop{\rm Tor}\nolimits}

\def\Ker{\mathop{\rm Ker}\nolimits}

\def\Im{\mathop{\rm Im}\nolimits}
\def\dim{\mathop{\rm dim_k}\nolimits}

\def\Fix{\mbox{{\rm Fix}}}

\begin{document}

\sf

\title{Cohomology of split algebras and of trivial extensions\thanks{This work
has been supported by the projects SECYT-ECOS A98E05 and SECYT-CAPES BR12/99/OG. The first author wishes to thank
FAPESP (Brazil) for financial support. The second author has a research scholarship from Cnpq (Brazil). The third
and fourth authors are research members of CONICET (Argentina).}}

\date{ }
\author{Claude Cibils,
 Eduardo Marcos,\\
 Mar\'\i a Julia Redondo  and   Andrea Solotar}

\maketitle

\begin{center}
{{\em{\rm To Idun Reiten for her 60th birthday}}}
\end{center}

\begin{abstract}
We consider associative algebras $\Lambda$ over a field provided with a direct sum decomposition of a two-sided
ideal $M$ and a sub-algebra $A$ -- examples are provided by trivial extensions or triangular type matrix algebras.
In this relative and split setting we describe a long exact sequence computing the Hochschild cohomology of
$\Lambda$. We study the connecting homomorphism using the cup-product and we infer several results, in particular
the first Hochschild cohomology group of a trivial extension never vanishes.

\end{abstract}

\small \noindent 2000 Mathematics Subject Classification : 16E40, 16D20

\noindent Keywords : cohomology, Hochschild, split, trivial extension.

\renewcommand\contentsname{\sf Contents}
\tableofcontents
\section {\sf Introduction }

In this paper we consider split algebras $\Lambda = A\oplus M$ where $A$ is a subalgebra of $\Lambda$ and $M$ is a
two--sided ideal. Our main purpose is to compute the Hochschild cohomology $H^*(\Lambda, \Lambda)$ using the
cohomology theory of $A$ and $M$.

Our motivations are three-fold. First, decompositions providing split algebras arise in various examples. Note for
instance that trivial extensions and triangular matrix algebras (see below) are split algebras, their cohomology
has been investigated recently by several authors \cite{begu,gms,mipl}. Second it is known that degree one
Hochschild cohomology provides insight to representation theory through universal covers, and its vanishing is
related to the notion of simply connected algebras, \cite{sk}. D. Happel shows in \cite{ha2} that a finite
representation type algebra over an algebraically closed field of characteristic zero is simply connected if and
only if the first Hochschild cohomology $H^1$ space of its Auslander algebra is zero. Moreover R. Buchweitz  and
S. Liu  provided a proof of the same statement assuming only that the field is algebraically closed (Oberwolfach
2000).
 It has been suspected that for a finite dimensional
algebra over an algebraically closed field, the vanishing of the first Hochschild cohomology space implies that
its ordinary quiver has no oriented cycles.  This was proved wrong and a family of counterexamples can be found in
\cite{buli}. It is also conjectured that a tilted algebra is simply connected if and only if its first Hochschild
cohomology space vanishes. This has been proved for tame tilted algebras, see \cite{asmape}. In another direction,
degree two cohomology concerns the deformation theory of algebras, see \cite{ge}. Finally split algebras are
interesting to study in relation to Happel's question \cite{ha}: if the Hochschild cohomology vector spaces of a
finite dimensional algebra vanish after some degree, is the algebra of finite homological dimension? The present
paper is a first step for considering this question in a relative and split framework.

We describe now the contents of each section of the article.

In section \ref{split} we obtain  a long exact sequence involving $H^*(\Lambda, \Lambda)$. If $M$ is a projective
left or right  $A$-module the other terms of this sequence are direct sums of vector spaces
$\Ext^q_{A-A}(M^{\otimes_A^p},X)$ for $p+q=*$ or $p+q=*+1$ depending on whether $X=M$ or $X=A$. We study in detail
the connecting homomorphism of this long exact sequence in order to obtain results on the dimensions of the
Hochschild cohomology vector spaces of $\Lambda$.

More precisely we obtain for any $\Lambda$--bimodule $X$ a double complex whose total cohomology is the Hochschild
cohomology of $\Lambda$ with coefficients in  $X$. The first quadrant spectral sequence involved converges, the
terms at the first level are unknown but can probably be approximated through new spectral sequences.
Nevertheless, if $M$ is projective as a left or right $A$--module, we show that the cohomology of the $p$--th
column ${\cal C}^p(X)$ is $\Ext^*_{A-A}(M^{\otimes_Ap},X)$.

In section \ref{connecting} we do not assume that $M$ is projective on one side. When $M$ is a square zero ideal
and the bimodule $X$ verifies $MX=XM=0$ then the horizontal differentials of the double complex are 0.
Consequently the Hochschild cohomology $H^*(\Lambda,X)$ is the direct sum of the cohomologies of the columns,
namely
$$H^n(\Lambda,X)= \bigoplus_{p+q=n}H^q({\cal
C}^p(X)).$$

Of course the $\Lambda$--bimodule $\Lambda$ does not verify the above hypothesis. Nevertheless the bimodules on
both sides of the sequence $ 0\to M \to \Lambda \to \Lambda/M \to 0$ do, and we consider the corresponding long
exact sequence in Hochschild cohomology. An interesting result we obtain is that the connecting homomorphism
$$\delta : H^n(\Lambda, \Lambda/M) \to
H^{n+1}(\Lambda, M)$$ is bigraded of bidegree $(1,0)$, that is $\delta = \oplus_{p+q=n} \delta^{p,q}$, where
$$\delta^{p,q}: H^q({\cal C}^p(\Lambda/M)) \to H^q({\cal C}^{p+1}(M)).$$

Moreover in section \ref{operations} we provide a precise description of $\delta^{p,q}$ involving only two terms,
using the cup product of a cocycle with the identity endomorphism of $M$. This description enables us to determine
whether $\delta^{p,q}$ annihilates or not in some interesting cases.

In section \ref{trivial} we consider a special case of split algebras: trivial extensions $TA$. They are of the
form $A\oplus DA$ where $DA$ is the dual $A$--bimodule  of $A$ endowed with the zero multiplicative structure. We
show that its first Hochschild cohomology is a direct sum of four vector spaces (see Theorem \ref{h1trivialext}).
One of the factors is  the center of $A$, in particular $H^1(TA, TA)$ never vanishes. More generally we show that
$H^n(A,A)\oplus H_n(A,A)$ is a direct summand of $H^n(TA,TA)$. This result is a consequence of the fact that for
these algebras the component $\delta^{0,q}$ of the connecting homomorphism is zero.

As a direct consequence of our computations we show that $H^1(TA,TA)=k \oplus H^1(A,A)$ for a \emph{one-way}
algebra $A$, see Definition \ref{oneway}. This generalizes previous results obtained in \cite{mj,ma}.

Finally we specialize to triangular matrix algebras and one--point extensions the results we have obtained for
general split algebras. In this way we recover computations performed in \cite{ci,gms,ha,mipl}.

\section {\sf Split algebras and the double complex}\label{split}

Let $k$ be a field. As stated in the introduction a split algebra $\Lambda$ is a $k$--algebra with a subalgebra
$A$ and a two--sided ideal $M$ such that $\Lambda=A \oplus M$. In other words $\Lambda$ consists of the following
data: a $k$--algebra $A$ and a multiplicative $A$--bimodule $M$ with a product, i.e. an associative $A$--bimodule
map $M \otimes_A M \to M$, $m \otimes m' \mapsto m.m'$. The algebra structure in $A\oplus M$ is given by
$$(a+m)(a'+m')=aa'+am'+ma'+m.m'.$$ Let $X$ be a
$\Lambda$--bimodule. As usual, the Hochschild cohomology vector spaces of $\Lambda$ with coefficients in $X$ are
the cohomology groups of the following  cochain complex, (see for instance \cite{ibra,caei,lo,we})
$$0\longrightarrow X \stackrel{d}{\longrightarrow} \Hom_k(\Lambda,X) \stackrel{d}{\longrightarrow} \cdots
\stackrel{d}{\longrightarrow} \Hom_k(\Lambda^{\otimes n},X) \stackrel{d}{\longrightarrow} \cdots $$ where for
$n\geq 1$
\begin{eqnarray*}
df(x_1 \otimes \dots \otimes x_{n+1}) &=& x_1 f(x_2\otimes \dots \otimes x_{n+1})\\ &+& \sum_{i=1}^n (-1)^i f(x_1
\otimes \dots \otimes x_ix_{i+1} \otimes \dots \otimes x_{n+1})\\ &+& (-1)^{n+1} f(x_1 \otimes \dots \otimes
x_{n})x_{n+1}
\end{eqnarray*}
and for $x\in X$ and $\lambda\in \Lambda$
$$(dx)(\lambda) = \lambda x - x\lambda .$$ Since
$\Lambda = A \oplus M$, we have a decomposition of $\Lambda^{\otimes n}$ as a direct sum of vector spaces in terms
of $A$ and $M$: let $M^{p,q}$ be the sub--vector space spanned by $(p+q)$--tensors $x_1 \otimes \dots \otimes
x_{p+q}$ such that exactly $p$ of the $x_i$'s belong to $M$ while the other $x_i$'s belong to $A$. Clearly
$$\Lambda^{\otimes n}=\bigoplus_{p+q=n}M^{p,q}.$$
Moreover the Hochschild complex above organizes in a double complex whose $(p,q)$-- spot is $\Hom_k(M^{p,q},X)$.
Indeed, the image of $d$ restricted to $\Hom_k(M^{p,q},X)$ is contained in $\Hom_k(M^{p+1,q},X) \oplus
\Hom_k(M^{p,q+1},X)$. The horizontal and vertical components of $d$ are denoted $d_h$ and $d_v$ respectively. The
cohomology of the total complex is $H^*(\Lambda,X)$.

\begin{prop}
\label{daigual} The vertical diffentials $d_v$ of the above double complex depend neither on  the product of $M$
nor on the actions of $M$ on $X$.
\end{prop}
\begin{Demo}
Let $f\in\Hom_k(M^{p,q},X)$, in other words $f: \Lambda^{\otimes (p+q)} \to X$ vanishes on $(p+q)$--tensors which
have not exactly $p$ components of $M$ and $q$ components of $A$. We evaluate $d_vf$  on a tensor $(x_1 \otimes
\dots \otimes x_{p+q+1})\in M^{p,q+1}$. We shall see that the terms where the product of $M$ or the action of $M$
on $X$ appear are zero. If $x_1\in M$ then $x_1f(x_2 \otimes \dots \otimes x_{p+q+1})=0$ since $(x_2 \otimes \dots
\otimes x_{p+q+1})\notin M^{p,q}$, regardless the action of $M$ on $X$. Similarly $f(x_1 \otimes \dots \otimes
x_ix_{i+1} \otimes \dots \otimes x_{p+q+1})=0$ if $x_i$ and $x_{i+1}$ belong to $M$ since the tensor belongs to
$M^{p-1,q+1}$ regardless the value of $x_ix_{i+1}$. The behaviour of the last term of the coboundary formula is
analogous to the first one. \qed
\end{Demo}\vskip5mm

In order to determine the vertical cohomology, we first note that the cohomology of the $0$--th column is the
Hochschild cohomology $H^*(A,X)=\Ext^*_{A-A}(A,X)$, where $X$ is considered as an $A$--bimodule by restriction of
scalars. From now on we simplify the notations omitting the tensor product sign for tensor products over the
ground field $k$; tensor signs between vectors are replaced by commas.

The following result is announced in \cite{ci}, but the proof provided there is incomplete.
\begin{teo}
\label{firstcolumn} The cohomology of the column $p=1$ is $\Ext^*_{A-A}(M,X)$.
\end{teo}

\begin{Demo}
We will first provide a free resolution of $M$ as an $A$--bimodule. Consider the bar resolution of $M$ as a left
$A$--module (see \cite[8.6.12]{we}),
$$\cdots\to AAM\to AM\to 0$$ and the Hochschild
resolution of $A$ as an $A$--bimodule $$\cdots\to AAA\to AA \to 0.$$ Tensoring them over $A$ provides the complex
$$\cdots \to AMAA\oplus AAMA \to AMA \to 0.$$ The cycles in each degree of the Hochschild resolution are
projective left $A$--modules since the resolution splits as a sequence of left $A$--modules. The K\"unneth formula
ensures that the last complex has zero homology in positive degrees and $M\otimes_A A = M $ in degree zero. Next
we apply the functor $\Hom_{A-A}(-,X)$ and we use the identification $\Hom_{A-A}(AZA,X)=\Hom_k(Z,X)$ in order to
verify that the coboundaries provide the first column of the double complex. \qed
\end{Demo}\vskip5mm

\begin{obser} \label{hom}
A direct computation shows that $\Hom_{A-A}(M^{\otimes_A p},X)$ is the zero degree cohomology of the $p$--th
column. In order to generalize this result to the other vertical cohomology groups for $p\geq 2$ we need to assume
additional hypothesis on $M$ as follows.
\end{obser}

\begin{teo} \label{ext}
Let $A$ be a $k$--algebra, $M$ be an $A$--bimodule which is right or left projective and $X$ be an $A\oplus
M$--bimodule. The cohomology of the $p$--th column in degree $q$ is $\Ext^q_{A-A}(M^{\otimes_Ap},X)$.
\end{teo}

The proof of the above theorem is given at the end of this section, as a consequence of the next result:

\begin{prop} \label{tor}
Let $C$, $B$ and $A$ be $k$--algebras, $_CN_B$ be a $C-B$--bimodule and $_BM_A$ be a $B-A$--bimodule. The homology
of the following complex is  equal to $\Tor^B_*(N,M)$
 $$ \cdots
\stackrel{b'}{\longrightarrow} CCNMA \oplus CNBMA \oplus CNMAA \stackrel{b'}{\longrightarrow} CNMA \longrightarrow
0,$$ where the term in degree $n$ is a free $C-A$--bimodule of the form $CZA$ and
$$Z=\bigoplus_{i+j+k=n;\ i,j,k \geq 0}
C^kNB^jMA^i.$$
\end{prop}

Note that the boundary formula is provided by the standard
 resolution of an algebra
as a bimodule:
\begin{eqnarray*}
b'(c,x_1, \dots, x_n, a) &=& (cx_1, \dots, x_n, a)\\ &+& \sum_{i=1}^{n-1}(-1)^{i} (c,x_1, \dots, x_ix_{i+1}, \dots
,x_n, a)\\ &+& (-1)^{n} (c,x_1, \dots, x_n a).
\end{eqnarray*}
By assumption products of type $nm$ are zero when $n \in N$ and $m \in M$; the bimodule action gives products in
case of elements of the form $bm$, $nb$, $cn$, etc. Each summand of the formula must have an element of $C$ and an
element of $A$ on each side, otherwise its value is $0$.

\begin{obser} Before proving  Proposition \ref{tor} we
note that for  a $C-A$--bimodule $X$, the functor $\Hom_{C-A}(-,X)$ applied to the above complex provides the
second column ${\cal C}^2(X)$ in case $A=B=C$ and $N=M$, since by virtue of Proposition \ref{daigual} we can
assume that the product of $M$ and the actions of $M$ on $X$ are trivial.
\end{obser}

\begin{Demo}{\sf (Proposition \ref{tor})}
We consider the bar resolution of $M$ as a left $B$--module $$ \cdots \longrightarrow BBM \longrightarrow BM
\longrightarrow 0$$ and we apply the functor $N\otimes_B-$

 $$ \cdots
\longrightarrow N \otimes_B BBM \longrightarrow N\otimes_B BM \longrightarrow 0,$$ obtaining in this way the
standard complex which is used to compute $\Tor^B_*(N,M)$
\begin{center}
\vskip -.5cm
\begin{equation}\label{cosa}
\cdots \longrightarrow NBM \longrightarrow NM \longrightarrow 0.
\end{equation}
\end{center}
Next we use the Hochschild resolution of $C$ as a $C$--bimodule $$ \cdots \longrightarrow CCC \longrightarrow CC
\longrightarrow 0.$$ Its homology is non--zero only in degree zero, with value $C$. Tensoring the above resolution
with (\ref{cosa}) over $C$ gives $$ \cdots \longrightarrow CNBM \oplus CCNM \longrightarrow CNM \longrightarrow
0.$$ We assert that the homology of this complex is still $\Tor^B_*(N,M)$. Indeed we can use again the K\"{u}nneth
formula since the set of cycles of the bar resolution of $C$ (which splits as a sequence of right $C$--modules) is
a projective right module, and the homology is zero except in degree zero with value $C$. Hence the homology of
the above complex is the tensor product of the homologies, that is $C\otimes_C\Tor^B_*(N,M)=\Tor^B_*(N,M)$.

Finally we consider the Hochschild resolution of $A$ and we tensor it over $A$ with the above complex. As before,
the resulting homology is $\Tor^B_*(N,M) \otimes_A A=\Tor^B_*(N,M)$. A non difficult computation shows that the
resulting boundaries coincide with those described in the statement.$\qed$
\end{Demo}\vskip5mm

\begin{prop} \label{tor2}
Let $A,B,C,D$ be $k$--algebras and $_DU_C$, $_CN_B$, $_BM_A$ be bimodules. Assume $\Tor^B_*(N,M)=0$ in positive
degrees. Then $\Tor^C_*(U, N\otimes_B M)$ is  the homology of the following complex: {\small $$ \cdots
\longrightarrow DDUNMA \oplus DUCNMA \oplus DUNBMA \oplus DUNMAA \longrightarrow DUNMA \longrightarrow 0$$} with
$n$--th term $DZA$ where
$$Z=\bigoplus_{i+j+k+l=n} D^lUC^kNB^jMA^i.$$

\end{prop}

\begin{Demo}
By hypothesis the complex of the preceding proposition has homology only in degree zero, with value $N\otimes_B
M$. Since the modules are $C$--free on the left, this complex is a resolution of the left $C$--module $N\otimes_B
M$. Applying the functor $U \otimes_C -$ provides a complex whose homology is $\Tor^C_*(U, N\otimes_B M)$.

We consider as in the previous proposition a resolution of the algebra $D$ as a $D$--bimodule. The K\"{u}nneth
formula shows that tensoring this resolution over $D$ with the complex obtained above provides a new complex whose
homology is $\Tor^C_*(U, N\otimes_B M)$.$\qed$
\end{Demo}\vskip5mm

\begin{obser}
If $\Tor^C_*(U, N\otimes_B M)$ is zero in positive degrees, the complex of Proposition \ref{tor2} becomes a
projective resolution of the $D-A$--bimodule $U\otimes_CN\otimes_BM$. Applying the functor $\Hom_{D-A}(-,_DX_A)$
to this projective resolution gives a cochain complex whose homology is $\Ext^*_{D-A}(U\otimes_C N \otimes_B
M,X)$.
\end{obser}

\begin{Demo}{\sf (Theorem \ref{ext})}
Since $M$ is right or left projective, then $\Tor^A_*(M,M)$ is zero in positive degrees. The complex of
Proposition \ref{tor}, in the case $N=M$, $A=B=C$, is a projective resolution of the $A$--bimodule $M\otimes_A M$.
We have already noticed that applying the functor $\Hom_{A-A}(-,X)$ to this resolution yields precisely the second
column of the double complex and, consequently, its cohomology is $\Ext^*_{A-A}(M\otimes_A M,X)$.

The same procedure applies to Proposition \ref{tor2}, and we obtain that the third column has cohomology
$\Ext^*_{A-A}(M^{\otimes_A 3},X)$. By induction the end of the proof is now obvious.$\qed$
\end{Demo}\vskip5mm

\section {\sf The connecting homomorphism}\label{connecting}

We return to the double complex which we use to compute the Hochschild cohomology of an arbitrary split algebra
$\Lambda = A \oplus M$ with coefficients in a $\Lambda$--bimodule $X$. The filtration of the total complex arising
from the columns provides a first quadrant spectral sequence, hence converging to $H^*(\Lambda, X)$ (see for
instance \cite{mcc,we}). In the preceding section we have computed the first level vector spaces $E^{*,*}_1$,
assuming $M$ is $A$--projective on one side. However the differential at the first level appears hard to compute
even with these hypothesis on $M$.

We focus on a special case of interest for specific computations that we will perform in the next section.

\begin{teo} \label{zero}
Let $A$ be a $k$-algebra, $M$ be an $A$-bimodule, $\Lambda = A \oplus M$ be the corresponding split algebra with
$M^2=0$, and let $X$ be a $\Lambda$--bimodule verifying $MX=XM=0$ (in other words $X$ is an $A$--bimodule with
actions trivially extended to $\Lambda$). Then the horizontal coboundaries of the double complex are zero.  As a
consequence
$$H^n(\Lambda, X)= \bigoplus_{p+q=n}H^q({\cal C}^p (X)),$$ where ${\cal C}^p(X)$ denotes the $p$-th column.
\end{teo}

\begin{Demo}
Let $\varphi : M^{p,q} \to X$ be a cochain, in other words $\varphi : \Lambda^{\otimes (p+q)} \to X$ is a cochain
that vanishes on each component $M^{p',q'} \neq M^{p,q}$ of $\Lambda^{\otimes (p+q)}$.

By definition $d_h \varphi = d \varphi \mid_{M^{p+1,q}}$ therefore
\begin{eqnarray*}
d_h \varphi (x_1, \dots , x_{p+q+1}) &=& x_1 \varphi (x_2, \dots , x_{p+q+1}) \\ &+& \sum_{i=1}^{p+q} (-1)^i
\varphi (x_1, \dots ,x_{i}x_{i+1}, \dots, x_{p+q+1})\\ &+& (-1)^{p+q+1} \varphi (x_1, \dots , x_{p+q})x_{p+q+1}.
\end{eqnarray*}
The first term is zero, indeed if $x_1 \in M$ we have $MX=0$, while if $x_1 \in A$ then $(x_2, \dots ,
x_{p+q+1})\in M^{p+1,q-1}$ and $\varphi$ is zero when evaluated on it. The last term is zero for the same
 reasons. Each middle term vanishes since either both $x_i$ and
$x_{i+1}$ belong to $M$ (hence $x_ix_{i+1}=0$) or $(x_1, \dots ,x_{i}x_{i+1}, \dots, x_{p+q+1}) \in M^{p+1,q-1}$.
$\qed$
\end{Demo}\vskip5mm

The above decomposition and the results of the previous section yield the following:

\begin{coro}
\label{suma} Let $A$ be a $k$--algebra and $M$ be an $A$--bimodule projective on one side. Let $\Lambda$ be the
split algebra $A \oplus M$ with $M^2=0$, and let $X$ be a $\Lambda$--bimodule such that $MX=XM=0$. Then
$$H^n(\Lambda,X)= \bigoplus_{p+q=n} \Ext^q_{A-A}(M^{\otimes_A p},X),$$ where $M^{\otimes_A 0}=A$.\end{coro}

We now consider for a split algebra $\Lambda = A \oplus M$ the following exact sequence of $\Lambda$--bimodules $$
0 \longrightarrow M \longrightarrow \Lambda \stackrel{\pi}{\longrightarrow} \Lambda/M=A \longrightarrow 0.$$
 Note that $M$ is a
$\Lambda$--bimodule since $M$ is a two--sided ideal of $\Lambda$. Of course $\Lambda/M$ is in fact $A$ considered
as a $\Lambda$--bimodule with zero actions of $M$ on both sides. This exact sequence of coefficients provides a
long exact sequence in Hochschild cohomology

\begin{eqnarray*}
0 \to & H^0(\Lambda, M)& \to H^0(\Lambda, \Lambda) \stackrel{\pi^0}{\to} H^0(\Lambda, A) \stackrel{\delta^0}{\to}
\\ & H^1(\Lambda,
M)& \to H^1(\Lambda, \Lambda) \stackrel{\pi^1}{\to} H^1(\Lambda, A) \stackrel{\delta^1}{\to}\\ & \dots &
\\ & H^n(\Lambda, M)& \to H^n(\Lambda, \Lambda)
\stackrel{\pi^n}{\to} H^n(\Lambda, A) \stackrel{\delta^n}{\to}
\\ &\dots
\end{eqnarray*}

Our next purpose is to describe the connecting homomorphism $\delta^n$ in order to combine this information with
knowledge on $H^*(\Lambda,A)$ and $H^*(\Lambda, M)$. This will provide information on $H^*(\Lambda, \Lambda)$
which is our main purpose. We begin by studying $\delta^0$.

\begin{prop}
\label{deltaescero} Let $A$ be a $k$-algebra, $M$ be an $A$-bimodule and let $\Lambda = A \oplus M$ be the
corresponding split algebra. The above connecting homomorphism $\delta^0$ vanishes if and only if the center $A^A$
of $A$ has symmetric action on M (i.e. $am=ma$ for every $a \in A^A$ and $m \in M$).
\end{prop}

\begin{Demo}
The center $\Lambda^{\Lambda}$ of $\Lambda=A \oplus M$ is as follows: $$\Lambda^\Lambda = [ A^A \cap A^M ] \
\oplus\ [M^M \cap M^A],$$ where $A^M$ are the elements of $A$ acting symmetrically on $M$, while $M^M$ is the
center of the multiplicative bimodule $M$ and $M^A= H^0(A,M)= \{m \in M \ \vert\  am=ma \ \mbox{for every $a \in
A$}\}$.

In the long exact sequence above $$0 \longrightarrow M^\Lambda \longrightarrow \Lambda^\Lambda
\stackrel{\pi^0}{\longrightarrow} A^\Lambda=A^A \stackrel{\delta^0}{\longrightarrow} H^1(\Lambda, M)
\longrightarrow \dots$$ we have $\Im \pi^0 = A^A \cap A^M$. Hence $\Ker \delta^0 = A^A \cap A^M$ so $\delta^0=0$
if and only if $A^{A}\cap A^M=A^{A}$ which is equivalent to $ A^A \subset A^M$.$\qed$
\end{Demo}\vskip5mm

\begin{ejem}
\label{deltacero} Let $f$ be an automorphism of $A$ and let $M=\ \! ^f\!\!A$ be the $A$--bimodule $A$ with left
action twisted by $f$. Then $\delta^0=0$ if and only if $f$ is the identity on central elements of $A$. Indeed,
let $\Fix f$ be the subalgebra of elements fixed by $f$. For $M=\ ^f\!\!A$ we have $A^M=A^A \cap \Fix f$.

\end{ejem}
\vskip1mm
\begin{ejem}
\label{dctrivial} In case $\Lambda$ is the trivial extension $TA$ of $A$ we have $M= DA$  and the center of $A$
acts symmetrically on $M$. Consequently  $\delta^0=0$.
\end{ejem}

We next prove that each connecting homomorphism has bidegree $(1,0)$.

\begin{prop}\label{diagonal}
Let $\Lambda =A \oplus M$ be a split algebra with $M^2=0$. The connecting homomorphism $$\delta^n : H^n(\Lambda,
\Lambda/M) \to H^{n+1}(\Lambda, M)$$ has bidegree $(1,0)$ with respect to the decomposition provided in Theorem
\ref{zero}.
\end{prop}

\begin{Demo}
We denote by $H^q({\cal C}^p(X))$ the cohomology of the $p$--th column in degree $q$, where $X$ is a
$\Lambda$--bimodule. Since we proved that the horizontal differentials are zero when $MX=XM=0$, we have for such
an $X$ $$ H^n(\Lambda, X) = \bigoplus_{p+q=n} H^q({\cal C}^p(X)).$$ Both $M$ and $\Lambda/M$ verify the above
assumption on $X$. Then $$\delta^n : \bigoplus_{p+q=n} H^q({\cal C}^p(\Lambda/M )) \to \left ( \bigoplus_{p+q=n}
H^q({\cal C}^{p+1}(M))\right ) \oplus H^{n+1}({\cal C}^0(M)).$$
 We assert that the image of
$\delta^n\mid_{H^q({\cal C}^p(\Lambda/M))}$ is contained in $H^q({\cal C}^{p+1}(M))$, hence $\delta^n=
\bigoplus_{p+q=n} \delta^{p,q}$ where
$$\delta^{p,q}: H^q({\cal C}^p(\Lambda/M)) \to
H^q({\cal C}^{p+1}(M)).$$ In order to prove the assertion let $\varphi: M^{p,q} \to A$ be a cocycle of the
Hochschild complex of $\Lambda/M$. We use the given inclusion $\Lambda/M \subset \Lambda=A \oplus M$ to obtain
$\overline \varphi: M^{p,q} \to \Lambda$, taking into account that $A$ inside $\Lambda$ has a non trivial action
of $M$ on it. The image of the coboundary of $\overline \varphi$ in the Hochschild complex of $\Lambda$ is
contained in $M$ and provides a well defined element in $H^{p+q+1}(\Lambda,M)$ by general arguments. Considering
the double complex for $\Lambda$, we have two components $d\overline\varphi = d_v\overline\varphi +d_h
\overline\varphi$. In fact $\delta \varphi = d\overline\varphi$. Now we will prove that $d \overline \varphi$ has
zero values on every component of $\Lambda^{\otimes n+1}$ except maybe on  $M^{p+1,q}$. In order to prove that
$d_v \overline\varphi=0$, let $(x_1, \dots , x_{p+q+1}) \in M^{p,q+1}$, then
\begin{eqnarray*}
d_v \overline\varphi (x_1, \dots , x_{p+q+1}) &=& x_1 \overline \varphi (x_2, \dots , x_{p+q+1}) \\ &+&
\sum_{i=1}^{p+q} (-1)^i \overline \varphi (x_1, \dots ,x_i x_{i+1}, \dots , x_{p+q+1}) \\ &+& (-1)^{p+q+1}
\overline \varphi (x_1, \dots , x_{p+q})x_{p+q+1}.
\end{eqnarray*}
We observe that the middle terms remain unchanged for $\varphi$ or $\overline\varphi$, namely {\small
$$\sum_{i=1}^{p+q} (-1)^i \overline \varphi (x_1, \dots ,x_i x_{i+1}, \dots , x_{p+q+1})\ = \ \sum_{i=1}^{p+q}
(-1)^i \varphi (x_1, \dots ,x_i x_{i+1}, \dots , x_{p+q+1}).$$} Concerning the first and the last terms we first
assume that both $x_1$ and $x_{p+q+1}$ belong to $A$, then all the terms of the sum are in $A$ and coincide with
the terms of $d_v\varphi$. Since $\varphi:M^{p,q}\to A$ is a cocycle, we obtain that the value of the above
expression is zero.

If $x_1\in M$ and $x_{p+q+1}\in A$ then $(x_2,\dots,x_{p+q+1})\in M^{p-1,q+1}$, hence $\varphi$ is zero  on it.
The last term remains in $A$, and all the terms of $d_v\overline\varphi$ evaluated on the tensor
$(x_1,\dots,x_{p+q+1})$ coincide with the terms of $d_v\varphi$ (the first one vanishes in both cases). Since
$\varphi$ is a cocycle, we infer that $d_v\overline\varphi(x_1,\dots,x_{p+q+1})=0$ also in this case. The
remaining cases $x_1,\ x_{p+q+1}\in M$, or $x_1\in A$, $x_{p+q+1}\in M$ can be studied in an analogous way.

We conclude that $\delta \varphi = d_h \overline\varphi \in \Hom(M^{p+1,q},M)$.$\qed$
\end{Demo}\vskip5mm

\section {\sf Operations} \label{operations}

We introduce operations in the double complex $\Hom_{A-A}(M^{*,*},X)$ of the previous sections in order to
describe the $(p,q)$--component $\delta^{p,q}$ of the connecting homomorphism
$$\delta^{p,q}: H^q({\cal C}^p(\Lambda/M)) \to
H^q({\cal C}^{p+1}(M)).$$ Recall that the following is an operation on Hochschild cohomology of bimodules over a
$k$--algebra $A$ (see \cite{caei,ge}). Let $X$ and $Y$ be $\Lambda$--bimodules, $f:\Lambda^{\otimes n} \to X$ and
$g:\Lambda^{\otimes m} \to Y$ be Hochschild cochains. The \emph{cup product} (see \cite{ge}) $f \smile g$ is
defined as the composition
$$\Lambda^{\otimes n+m} \cong \Lambda^{\otimes n}
\otimes \Lambda^{\otimes m} \stackrel{f\otimes g}{\longrightarrow}  X \otimes Y \longrightarrow X \otimes_\Lambda
Y .$$ One has $d(f \smile g)= df \smile g + (-1)^n f \smile dg \ $, so the product is well defined in cohomology:
$$H^n(\Lambda,X) \otimes H^m(\Lambda,Y) \to H^{n+m}(\Lambda,X\otimes_\Lambda Y).$$

In case $\Lambda=A \oplus M$ is a split algebra this operation goes clearly through the double complex, that is,
if $f:M^{p,q} \to X$ and $g:M^{p',q'} \to Y$ are cochains, then $$f \smile g: M^{p,q}\otimes M^{p',q'} \to X
\otimes_\Lambda Y$$ is the product cochain. Note that $M^{p,q}\otimes M^{p',q'}$ is naturally a direct summand of
$M^{p+p',q+q'}$, the value of $f\smile g$ on the complement is zero.

\begin{teo}
\label{smile} Let $\Lambda=A\oplus M$ be a split algebra with $M^2=0$, and let $$0 \to M \to \Lambda \to \Lambda/M
\to 0$$ be the corresponding short exact sequence. The $(p,q)$--component $\delta^{p,q}$ of the connecting
homomorphism $\delta$ in the long exact Hochschild cohomology sequence of $\Lambda$ is given by $$\delta^{p,q}
\varphi\ =\ 1_M \smile \varphi\ \ +\ \ (-1)^{p+q+1} \varphi \smile 1_M.$$
\end{teo}

\begin{obser}
In the statement of this theorem, $\varphi : M^{p,q} \to \Lambda/M$ is an arbitrary cocycle and $1_M : M \to M$ is
the identity morphism which is indeed a $1$--cocycle; it belongs to the $(1,0)$--spot of the double complex and
corresponds to the projection $A\oplus M\to M$ in the usual Hochschild complex of $M$. Note that if $M\neq 0$ this
projection is a non--zero element in $H^1(\Lambda, M)$.

Note also that we have $\Lambda/M \otimes_\Lambda M = \Lambda/M \otimes_A M = A\otimes_A M = M$ as well as $M
\otimes_\Lambda \Lambda/M =M$.
\end{obser}

\begin{Demo}
We lift the cocycle $\varphi$ to $\overline \varphi : M^{p,q} \to \Lambda$ as in the previous section. Since
$\delta^{p,q} \varphi = d_h \overline \varphi$, we consider $(x_1, \dots, x_{p+q+1}) \in M^{p+1,q}$. In the
coboundary formula the middle terms are all zero. If $x_i$ and $x_{i+1}$ belong to $M$ then $x_ix_{i+1}=0$ since
$M^2=0$. Otherwise $x_i$ or $x_{i+1}$, or both of them lie in $A$, hence $(x_1, \dots,x_ix_{i+1}, \dots,
x_{p+q+1})$ belongs to $M^{p+1,q-1}$ and $\varphi$ is zero evaluated on this tensor. We have proved that
\begin{eqnarray*}
\delta^{p,q} \varphi (x_1, \dots, x_{p+q+1}) &=& x_1 \varphi (x_2, \dots, x_{p+q+1}) \\ &+& (-1)^{p+q+1} \varphi
(x_1, \dots, x_{p+q})x_{p+q+1},
\end{eqnarray*}
which corresponds to the formula involving the cup product with the identity endomorphism of $M$. $\qed$
\end{Demo}\vskip5mm

\begin{ejem}
We describe the connecting homomorphism component
$$\delta^{p,0} : \Hom_{A-A}(M^{\otimes_Ap},A) \to
\Hom_{A-A}(M^{\otimes_A p+1}, M).$$ Let $\varphi \in \Hom_{A-A}(M^{\otimes_Ap},A) $ be a cocycle, then
$$(\delta^{p,0}\varphi) (m_1, \dots, m_{p+1})= m_1 \varphi (m_2, \dots, m_{p+1})+ (-1)^{p+1} \varphi (m_1, \dots,
m_{p})m_{p+1}.$$ For $p=0$ we have $$\delta^{0,0} : \Hom_{A-A}(A,A) \to \Hom_{A-A}(M, M).$$ Recall that
$A\otimes_A M$ and $M \otimes_A A$ are identified with $M$. Then $\delta^{0,0}\varphi(m)= m \varphi(1)-\varphi(1)
m$. Since the center $A^A$ is identified with $\Hom_{A-A}(A,A)$, then $\Ker \delta^{0,0}=A^A \cap A^M$, in other
words the kernel of $\delta^0$ is the set of central elements of $A$ which act symmetrically on $M$, as in
Proposition \ref{deltaescero}. Note also that $\delta^0=\delta^{0,0}$.
\end{ejem}

\section {\sf Trivial extensions}\label{trivial}

\begin{defi}
The \emph{trivial extension} $TA$ of an algebra $A$ is the split algebra obtained by using the $A$--bimodule
$DA=\Hom_k(A,k)$ endowed with the zero multiplicative structure.
\end{defi}

We recall that for trivial extensions the connecting homomorphism $\delta^0$ of the long exact cohomology sequence
is zero, see Example \ref{dctrivial}. Our next purpose is to compute the first Hochschild cohomology vector space
of a trivial extension.  For this we study the first connecting homomorphism $\delta^1$. Since $\delta^0=0$ the
long exact cohomology sequence for $TA$ gives the following exact sequence$$
\begin{array}{llllllll}
& 0\rightarrow &H^1(TA,DA) & \rightarrow &H^1(TA,TA) & \rightarrow & H^1(TA, A) &\stackrel{\delta^1}{\rightarrow}
\\&& H^2(TA,DA)& \rightarrow &H^2(TA,TA) & {\rightarrow}
&H^2(TA,A) &\stackrel{\delta^2}{\rightarrow}\\ && \dots $$
\end{array}$$
Using Theorem \ref{zero} for $X=TA/DA=A$ or $X=DA$ we have $$H^n(TA, X)= \bigoplus_{p+q=n}H^q({\cal C}^p (X)).$$
Moreover $\delta^n=\oplus_{p+q=n}\delta^{p,q}$ where
$$\delta^{p,q}:H^q({\cal C}^p (A)) \longrightarrow
H^q({\cal C}^{p+1} (DA)) $$

\begin{obser}\label{bajosgrados}
The following facts hold without any projectivity hypothesis on the $A$--bimodule $DA$

\begin{enumerate}
\item  $H^*({\cal C}^0(X))=H^*(A,X)$,
\item  $H^*({\cal C}^1(X))=\Ext^*_{A-A}(DA,X)$ (cf. Theorem
\ref{firstcolumn}),
\item  $H^0({\cal C}^p(X))=\Hom_{A-A}(DA^{\otimes_Ap},X)$ (cf. Remark \ref{hom}).
\end{enumerate}

 Using them we get that for $n=1$
  {\small

$
\begin{array}{ll}
\delta^1:&\Hom_{A-A}(DA,A)\ \oplus\   H^1(A,A) \stackrel {\mbox{\tiny$\left(\begin{array}{cc}\delta^{1,0} & 0\\
0 &\delta^{0,1}\\0 & 0\end{array}\right)$}} {\ \vector(1,0){80}}\\ &\Hom_{A-A}\left(DA\otimes_ADA,DA\right)\
\oplus\ \Ext^1_{A-A}(DA,DA)\ \oplus\ H^2(A,DA)
\end{array}
$}

\end{obser}

\begin{prop}\label{deltas}
The connecting morphisms  $\delta^{0,1}$ and $\delta^{1,0}$ verify:
\begin{itemize}
\item[i)] $\delta^{0,1}=0$
\item[ii)] Under appropriate identifications $\delta^{1,0}\varphi=\varphi+\varphi^*$
\end{itemize}
\end{prop}

\begin{obser}The first item of this proposition will be generalized
in Proposition \ref{q}.
\end{obser}

\begin{Demo}
\begin{itemize}
\item[i)]Let $\varphi$ be a vertical cocycle at the $(0,1)$--spot of
the double complex of $A$, namely $\varphi :A\to A$ is a usual derivation of the algebra $A$. We know that
$\delta^{0,1}\varphi=1_{DA}\smile\varphi\ +\ \varphi\smile 1_{DA}$ (see Theorem \ref{smile}):
$$\begin{array}{cccccl}\delta^{0,1}\varphi :
&(DA)A\ &\oplus \ &A(DA) &\longrightarrow &DA\\ &(f,a)&+&(b,g)&\mapsto &f\varphi(a)+\varphi(b)g.
\end{array}$$
We assert that $\delta^{0,1}\varphi$ is actually a vertical coboundary in the double complex of $DA$, namely
$\delta^{0,1}\varphi =d_v\varphi^*$. Indeed
$$\left( d_v\varphi^*\right)(f,a)=
-\varphi^*(fa)+\varphi^*(f)a.$$ For every $x\in A$ we have $$
\begin{array}{rllll}
-\left( \varphi^*(fa) \right)(x)\ &+\  & \left( \varphi^*(f)a \right)(x)&=&\\ -(fa)\left( \varphi(x) \right) &+
&\varphi^*(f)(ax) &=&\\ -f\left( a\varphi(x) \right) &+ &f(\varphi(ax)) &=&\\ f\left( -a\varphi(x)\ \right. & +
&\left.  \ \varphi(ax)\right) &=& f(\varphi(a)x).
\end{array}
$$ The last equality holds since $\varphi$ is a
derivation. Finally we obtain
$$\left(d_v\varphi^*\right)(f,a)=f\varphi(a).$$
Similarly we prove that $\left(d_v\varphi^*\right)(b,g)$ equals $\varphi(b)g$.

\item[ii)]
Let $\varphi \in \Hom_{A-A}(DA,A)$ or by adjointness let $\beta:DA\otimes_{A-A}DA \to k$ be a bilinear form, given
by $\beta(f,g)=g\left(\varphi (f)\right)$.  We know that $\delta^{1,0}\varphi=1_{DA} \smile \varphi + \varphi
\smile 1_{DA}$, more precisely $(\delta^{1,0}\varphi)(f,g)=f\varphi(g)+\varphi(f)g$. Now each $\psi \in
\Hom_{A-A}(DA\otimes_ADA,DA)$ is also identified with a bilinear form $\beta:DA\otimes_{A-A}DA \to k$, namely
$\beta(f,g)=\psi(f,g)(1)$.  Through this identification, we have that $\delta^{1,0}\beta=\beta+\beta^t$, where
$\beta^t(f,g)=\beta(g,f)$.  Indeed $$(\delta^{1,0} \varphi )(f,g)(1)=\left( f \varphi (g) \right )(1)+\left(
\varphi(f)g \right)(1)=f(\varphi (g))+g(\varphi(f)).$$
\end{itemize}\qed
\end{Demo} \vskip5mm

We consider the set  of skew--symmetric bilinear forms $\beta$ over $DA$ such that $\beta(fa,g)=\beta(f,ag)$ and
we denote this set $\mathrm{Alt}_A(DA)$. In the proof of the next Theorem we will show that
$\mathrm{Alt}_A(DA)=\mathrm{Ker}\delta^{1,0}$. This vector space coincides with $\mathcal{E}(DA)$ as considered by
M. Saorin in \cite{ma}.

We use a star symbol in order to denote the dual of a vector space, while the notation $D$ is kept for the dual of
a vector space endowed with a bimodule structure.

\begin{teo}\label{h1trivialext}
Let $TA$ be the trivial extension of a finite--dimensional algebra $A$. Then
$$H^1(TA,TA)\ \ =\ \ A^A\ \oplus\ H_1(A,A)^*\
\oplus \ H^1(A,A)\ \oplus\ \mathrm{Alt}_A(DA).$$
\end{teo}

Before proving this result we note that since the center of a $k$-algebra is not zero we get the following result.

\begin{coro}
Let $TA$ be the trivial extension of a finite--dimensional algebra $A$. Then the first Hochschild cohomology group
of $TA$ do not vanish.

\end{coro}

\begin{Demo}{\sf Theorem \ref{h1trivialext}.}
From the long exact sequence and the description of $\delta^1$, we have $$ H^1(TA,TA)\ = \ H^1(TA,DA) \oplus
\mathrm{Ker}\delta^1.$$ We have that $\mathrm{Ker}\delta^1 =H^1(A,A) \oplus \mathrm{Ker} \delta^{1,0}$ and
$$\mathrm{Ker} \delta^{1,0}= \left\{\varphi\in\Hom_{A-A}(DA,A)\ \mid \varphi+\varphi^*=0 \right\}.$$ Using
adjointness we have $$\Hom_{A-A}(DA,A)=\left(DA\otimes_{A-A} DA\right)^*.$$ So we have proved the following
$$\mathrm{Ker}\delta^{1,0}\ =\ \mathrm{Alt}_A(DA).$$ From Remark \ref{bajosgrados} $$H^1(TA,DA)\ =\
\Hom_{A-A}(DA,DA)\ \oplus\ \Ext^1_{A-A}(A,DA).$$ Actually
$$\Hom_{A-A}(DA,DA)=\Hom_{A-A}(A,A)=A^A.$$ Finally
$$\ H^1(A,DA)\ =\ H_1(A,A)^*$$ since for a finite
dimensional algebra $A$ and a finite dimensional $A$--bimodule $N$ the following fact holds: $H_n(A,N)^*\ =\
H^n(A,DN)$ (see for instance \cite{caei}). Note also that $N$ and $DDN$ are bimodules which are canonically
isomorphic by the evaluation map. \qed
\end{Demo}\vskip5mm

\begin{teo}
Let $A$ be an arbitrary algebra (not necessarily finite--dimensional). Then $$H^1(TA,TA)\ \ =\ \ A^A\ \oplus\
H^1(A,TA)\ \oplus\ \mathrm{Alt}_A(DA).$$
\end{teo}

\begin{Demo}
In the proof of the above theorem note that $H^1(A,A)\oplus H^1(A,DA)=H^1(A,TA)$. \qed
\end{Demo} \vskip5mm

\begin{teo}\label{sumando} Let $A$ be a finite dimensional $k$--algebra.
Then the vector space $H^n(A,A)\oplus H_n(A,A)$ is a direct summand of $H^n(TA,TA)$.
\end{teo}

In order to prove this theorem, we provide the following result generalizing Proposition \ref{deltas}. Recall that
$\delta^{0,q}: H^q(A,A) \to \Ext^q_{A-A}(DA,DA)$ is a component of the connecting homomorphism $\delta^q:
H^q(TA,A) \to H^{q+1}(TA,DA)$.

\begin{prop}\label{q}
We have $\delta^{0,q}=0$ for all $q\geq 1$.
\end{prop}

\begin{Demo}
Let $\phi \in \Hom(A^{\otimes q},A)$ be a Hochschild cocycle at the $(0,q)$--spot of the double complex.  Next we
provide $\phi' \in \Hom((DA)^{1,q-1},DA)$ such that $d_v\phi'=\delta^{0,q}\phi$, by the following formula
$$\phi'(a_1, \dots, a_n, f, b_1, \dots, b_m)(x)= \epsilon(n,q)f \left (\phi( b_1, \dots, b_m, x,  a_1, \dots,
a_n)\right)$$ where $n+m+1=q$ and
\[ \epsilon(n,q)= \left\{
\begin{array}{ll}
-1 \qquad & \mbox{if $n$ is odd}\\ (-1)^{q+1} & \mbox{if $n$ is even.}
\end{array} \right. \]
Observe that $\delta^{0,q}(\phi) \in \Hom((DA)^{1,q-1},DA)$ and $\delta^{0,q}(\phi)=1 \smile \phi + (-1)^{q+1}
\phi \smile 1$. So the following three cases arise: {\small \[
\begin{array}{ll}
&(\delta^{0,q}\phi )(f, b_1,  \dots, b_q)(x)=f\left( \phi(b_1, \dots, b_q)x\right)  \\
\\ &(\delta^{0,q}\phi )(a_1, \dots, a_i ,f, b_1,
\dots, b_j)(x)=0 \qquad \mbox{for $i \neq 0$, $j \neq 0$, $i+j=q$}\\ \\ &(\delta^{0,q}\phi )(a_1, \dots,
a_q,f)(x)=(-1)^{q+1} f\left( x\phi(a_1, \dots, a_q)\right).
\end{array}\]}
The verification that $d_v\phi'=\delta^{0,q}\phi$ is left to the reader.$\qed$
\end{Demo}\vskip5mm

\begin{Demo} (\sf Theorem \ref{sumando})
By the previous result $H^q(A,A)$ is contained in the image of the morphism $H^q(TA,TA)\to H^q(TA,A)$. Concerning
the homology factor note, as we remarked before, that $H^q({\cal C}^0(DA))=H^q(A,DA)$. We know from Proposition
\ref{diagonal} that $H^q(A,DA) \cap \Im \delta^{q-1}=0$ therefore $H_q(A,A)=H^q(A,DA) \subset H^q(TA,TA)$. $\qed$
\end{Demo}\vskip5mm

We generalize now a result obtained in \cite{mj,ma}. In \cite{mj} it is shown that for a triangular schurian
algebra $A$ we have $H^1(TA,TA)=k \oplus H^1(A,A)$. Actually the same equality holds for triangular algebras, or
for $2$-nilpotent algebras whose quiver do not contain oriented cycles of length $\leq 2$, see \cite{ma}.

Our next purpose is to use our previous computations on $H^1$ of a trivial extension in order to show that this
result holds  for one-way algebras, a family of algebras that we define below and which include the algebras
considered above. Note that the proof in \cite{ma} of the above equality under the mentioned hypothesis
 also works for one-way algebras.

\begin{defi}\label{oneway}
A {\em one-way} algebra is a finite dimensional algebra endowed with a complete set $S$ of orthogonal idempotents
such that
\begin{enumerate}
\item
For $e\neq f$ in $S$, if $eAf \neq 0$ then $fAe=0$.
\item
For all $e\in S$, we have $\dim(eAe)=1$.
\item
$S$ has more than one element (\emph{ie.} $A$ is not k) and $A$ is an indecomposable algebra (\emph{ie.} the graph
with set of vertices $S$ and an edge between $e$ and $f$ in case $eAf$ or $fAe$ is not zero is a connected graph).
\end{enumerate}
\end{defi}

\begin{teo}\label{one-way}
Let $A$ be a finite dimensional one-way algebra, and let $TA$ be its trivial extension. Then
$$H^1(TA,TA)\ \ =\ \ k \ \oplus\ H^1(A,A)\ .$$
\end{teo}

In order to prove this formula we use Theorem \ref{h1trivialext}, and two  results as follows.

\begin{lem}\label{homDA,A}
Let $A$ be a one-way algebra. Then
$$\Hom_{A-A}(DA,A)=0.$$
\end{lem}
\begin{Demo}
Take $\varphi \in \Hom_{A-A}(DA, A)$ and $e, f \in S$ distinct.  Since $\varphi (D(eAf)) \subset fAe$, the form
$\varphi$ has to vanish on $D(eAf)$. Now $\varphi(D(eAe)) \subset eAe$ so $\Im \varphi \subset \oplus_{e \in S}
eAe$. But since the algebra is indecomposable and different from $k$ there exist $f \neq e$ such that $eAf \neq
0$. Hence there is no non--zero two-sided ideal contained in the vector space $\oplus_{e\in S} eAe$. \qed
\end{Demo} \vskip5mm

\begin{lem}
Let $A$ be a one-way algebra for a system $S$ of idempotents. Then $H_1(A,A)=0 .$
\end{lem}
\begin{Demo}
Let $E=\times_{e\in S}ke$ be the subalgebra of $A$ generated by $S$. Note that $A\otimes_E A$ is a projective
$A$-bimodule since $$ A\otimes_E A = \bigoplus_{e\in S}Ae\otimes eA$$  and each summand is a projective
$A$-bimodule using the fact that
$$A\otimes A= \bigoplus_{e,f\in S}Ae\otimes fA.$$

Note also that the  decomposition $$A\otimes_E A\otimes_E A=\bigoplus_{e,f \in S}Af\otimes fAe\otimes eA$$ shows
that this bimodule is projective as an $A$-bimodule.

Consider now the projective resolution of $A$ as an A-bimodule $$ \dots \to A\otimes_E A\otimes_E A\to A\otimes_E
A \to A \to 0$$ where the boundary formula is provided by the standard Hochschild resolution of $A$ as an
$A$--bimodule. The homotopy contraction showing the exactness is defined as usual, by inserting $1$ at the
beginning of each tensor.

Applying the functor $ - \otimes_{A^e}A$ one gets, after decomposing in terms of the orthogonal idempotents
$$\dots \to \bigoplus_{e,f,g\in S}eAf\otimes fAg\otimes gAe  \to \bigoplus_{e,f\in S}eAf\otimes fAe \to
\bigoplus_{e\in S}eAe \to 0.$$

But $\bigoplus_{e,f\in S}eAf\otimes fAe = \bigoplus_{e\in S}eAe\otimes eAe $, since $eAf\neq 0$ implies $fAe=0$ if
$f\neq e$. Also, if $e$ is any primitive idempotent, $eAe$ is isomorphic to $k$. The boundary map, using this
isomorphism, is null on $\bigoplus_{e\in S}eAe\otimes eAe$.

As before, the following term $\bigoplus_{e,f,g\in S}eAf\otimes fAg\otimes gAe$ of the complex may be written as:
$$\left(\bigoplus_{e\in S}eAe\otimes eAe \otimes
eAe\right)\oplus \left(\bigoplus_{e\neq f\in S}eAf\otimes fAe \otimes eAe\right)\oplus$$ $$ \left(\bigoplus_{e\neq
f\in S}eAe\otimes eAf \otimes fAe\right)\oplus \left(\bigoplus_{e\neq f \neq g\in S}eAf\otimes fAg \otimes
gAe\right).$$ The second and third summands are zero, and the restriction of the boundary map to the first one,
composed with the isomorphism $eAe\cong k$ is the identity. So, already restricted to this first summand, the
boundary map is surjective. Then $H_1(A,A)=0$.
 \qed
\end{Demo}
\vskip5mm
\begin{Demo}(\sf Theorem \ref{one-way})
We recall the decomposition of Theorem \ref{h1trivialext}: $$H^1(TA,TA)\ \ =\ \ A^A\ \oplus\ H_1(A,A)^*\ \oplus \
H^1(A,A)\ \oplus\ \mathrm{Alt}_A(DA).$$ The hypothesis on $A$ implies that the center $A^A$ of $A$ is the field
$k$. From  Lemma \ref{homDA,A} we get $\mathrm{Alt}_A(DA)=0$, since $\mathrm{Alt}_A(DA)\subset \Hom_{A-A}(DA,A)$.
The previous theorem shows that $H_1(A,A)=0$. \qed \vskip5mm

\end{Demo}

The aim of the last part of this section is to show that the connecting homomorphisms of the long exact sequence
on Hochschild cohomology are not all zero in general.

 We consider split algebras with $M=A$ and $M^2=0$. These
algebras are isomorphic to $A[x]/<x^2> \simeq A \otimes k[\epsilon]$, where $k[\epsilon]=k[x]/<x^2>$ is the
algebra of dual numbers. We denote them by $A[\epsilon]$. Recall that an algebra $A$ is {\it symmetric} if $A$ is
isomorphic to $DA$ as an $A$--bimodule. In this case the trivial extension $TA$ of $A$ coincides with the split
algebra $A[\epsilon]$.

It is well known (see for instance \cite{caei}, \cite[Prop. 9.4.1]{we}) that if $A$ and $B$ are $k$--algebras (one
of them finite dimensional) we have $$H^n(A \otimes B, A \otimes B)= \bigoplus_{p+q=n}H^p(A,A) \otimes H^q(B,B).$$
It is also well known that if $k$ is of characteristic different from $2$ then
\[ \dim H^*(k[\epsilon], k[\epsilon])= \left\{
\begin{array}{ll}
2 \quad & \mbox{if $*=0$} \\ 1 & \mbox{if $* > 0$.}
\end{array} \right. \]
If $\mathrm{char}\ k = 2$ then $\dim H^n(k[\epsilon], k[\epsilon])=2$ for all $n$. For a $k$--algebra $A$ we infer
that in characteristic different from $2$
$$H^n(A[\epsilon],A[\epsilon])=H^n(A) \oplus
\left( \bigoplus_{i=0}^n H^i(A) \right),$$ while in characteristic $2$
$$H^n(A[\epsilon],A[\epsilon])= \bigoplus_{i=0}^n
\left( H^i(A) \oplus H^i(A)\right).$$

Let $\Lambda=A \oplus M$ be a split algebra with $M^2=0$. Assume that all the connecting homomorphisms are zero.
Then
$$H^n(\Lambda,\Lambda)= \left(\bigoplus_{p+q=n}
H^q({\cal C}^p(M))\right) \oplus \left(\bigoplus_{p+q=n} H^q({\cal C}^p(A))\right).$$ In case $M$ is projective on
one side and all connecting homomorphisms are zero, we get $$H^n(\Lambda,\Lambda)= \left(\bigoplus_{p+q=n}
\Ext_{A-A}^q(M^{\otimes_Ap},M) \right) \oplus \left(\bigoplus_{p+q=n} \Ext_{A-A}^q(M^{\otimes_Ap},A)\right).$$ In
case $M=A$, and still assuming that all connecting homomorphisms are zero, we get
$$H^n(A[\epsilon],A[\epsilon])= \bigoplus_{i=0}^n
\left( H^i(A) \oplus H^i(A)\right)$$ which holds only in characteristic two. Hence the connecting homomorphisms
are not zero in general.

\begin{obser}
For trivial extensions one can describe the component $$\delta^{p,0}: \Hom_{A-A}(DA^{\otimes_Ap},A)\rightarrow
\Hom_{A-A}(DA^{\otimes_Ap+1},DA)$$ of the connecting homomorphism as follows, generalizing the second item of
Proposition \ref{deltas}. The cyclic group of order $p+1$ acts on $\Hom_{A-A}(DA^{\otimes_Ap+1},DA)$ {\em via}
$$(t\varphi)(f_1\otimes\cdots\otimes
f_{p+1})=\varphi(f_2\otimes\cdots f_{p+1}\otimes f_1).$$ Identifying by adjunction the source with the target of
$\delta^{p,0}$ we obtain
$$\delta^{p,0}\varphi = t\varphi +
(-1)^{p+1}\varphi.$$

\end{obser}

\section{\sf Triangular matrix algebras and one--point extensions}

Recall that a \emph{triangular matrix algebra} {\small $\left(
\begin{array}{cc}
A & 0\\ M & B
\end{array}\right)$}
consists of two algebras $A$ and $B$ and a $B-A$ bimodule $M$, the product is obtained by matrix multiplication.
Note that in case $B$ is the ground field $k$ such algebras are called \emph{one--point extensions} of $A$. Our
next purpose is to specialize to these algebras the results we have obtained for split algebras in order to
recover results of C. Cibils, S. Michelena and M.I. Platzeck in \cite{ci,mipl}, and by D. Happel for one--point
extensions \cite{ha}, see also \cite{begu,gms}.

\begin{obser}
Triangular matrix algebras are split algebras with zero bimodule product. Indeed consider the algebra $A\times B$
and the trivially extended $A\times B$--bimodule $M$ with structure given by $(a,b)m=bm$ and $m(a,b)=ma$. The
split algebra $(A\times B)\ \oplus\ M$ with $M^2=0$ is exactly the algebra \small $\left(
\begin{array}{cc}
A & 0\\ M & B
\end{array}\right).$
\normalsize
\end{obser}

Let $T=$ {\small $\left(
\begin{array}{cc}
A & 0\\ M & B
\end{array}\right)$}
 be a triangular matrix algebra. Consider the exact sequence of $T$--bimodules
 $$0\to M\to T\to A\times B\to 0$$ and the corresponding long exact sequence in Hochschild cohomology

$$
\begin{array}{lclclclcl}
0\rightarrow &H^0(T,M) & \rightarrow &H^0(T,T) & \rightarrow & H^0(T, A\times B)&\stackrel{\delta^0}{\rightarrow}
\\& H^1(T,M)& \rightarrow &\cdots \\
&&&\cdots &\rightarrow&H^{n-1}(T, A\times
B)&\stackrel{\delta^{n-1}}{\rightarrow}\\&H^n(T,M)&\rightarrow&H^n(T,T)&\rightarrow&H^n(T,A\times
B)&\stackrel{\delta^{n}}{\rightarrow}\\ &H^{n+1}(T,M)&\rightarrow&\cdots

\end{array}$$

We will use a suitable version of Corollary \ref{suma} in order to describe $H^n(T,M)$ and $H^n(T,A\times B)$. The
following fact will enable us to perform a $\Tor$ computation for recovering Cibils and Michelena-Platzec Theorem.

\begin{lem}
If $M$ is projective as a left $B$--module, the trivially extended $A\times B$--bimodule $M$ is a projective left
$A\times B$--module.
\end{lem}
\begin{Demo}
Note that $B$ is projective as a left $A\times B$--module, consequently the same  holds for a direct summand of a
free $B$--module.\hfill\qed
\end{Demo}

\bigskip

The next result simplifies considerably this description.

\begin{lem}\label{tensorzero}
Let $M$ be a $B-A$--bimodule trivially extended to an $A\times B$--bimodule. Then for $p\geq 2$ we have
$$M^{\otimes^p_{(A\times B)}}=0.$$
\end{lem}
\begin{Demo}
For $m\in M$ and $n\in M$ we have $$m\otimes n \ =\ m(1,0)\otimes n \ =\  m\otimes (1,0)n \ =\ m\otimes 0\ =\
0.$$\hfill\qed
\end{Demo}

\begin{teo}
(see \cite{ci,mipl}) Let $A$ and $B$ be $k$--algebras, $M$ be a $B-A$--bimodule and let $T=$\footnotesize{ $\left(
\begin{array}{cc}
A & 0\\ M & B
\end{array}\right)$}
 \normalsize
 $= (A\times B) \oplus M$ be the
triangular matrix algebra or equivalently the corresponding split
 algebra. Then there is a long exact sequence in Hochschild
 cohomology
\small{
 $$
\begin{array}{lcclllll}
0\rightarrow &0 & \rightarrow &H^0(T,T) & \rightarrow & H^0(A, A)\oplus H^0(B,B)
&\stackrel{\delta^{0,0}}{\rightarrow}

\\& \Hom_{B-A}(M,M)& \rightarrow &\cdots \\

&&&\cdots &\rightarrow &H^{n-1}(A,A)\oplus H^{n-1}(B,B) &\stackrel{\delta^{n-1,0}}{\rightarrow}\\

&\Ext^{n-1}_{B-A}(M,M)&\rightarrow&H^n(T,T)&\rightarrow&H^n(A,A)\oplus H^n(B,B)
&\stackrel{\delta^{n,0}}{\rightarrow}\\ &\Ext^{n}_{B-A}(M,M)&\rightarrow&\cdots

\end{array}$$
}
\end{teo}
\begin{Demo}
Theorem \ref{zero} provides the following decompositions

$$H^n(T,A\times
B)=\bigoplus_{p+q=n}H^q({\cal C}^p (A\times B))$$

$$H^{n+1}(T,M)=H^{n+1}(A\times B,M)\ \oplus\ \bigoplus_{p+q=n+1}H^q({\cal C}^{p+1} (M))
$$ and Proposition \ref{diagonal} shows that the connecting homomorphism $\delta^n$ is bigraded of
bidegree $(1,0)$, that is, $\delta^n=\bigoplus_{p+q=n}\delta^{p,q}.$

In Section \ref{split} we have proved that the cohomology in degree $q$ of the column ${\cal C}^p (X)$ for a
$T$-bimodule $X$ is $\Ext^q_{(A\times B)-(A\times B)}\left(M^{\otimes^p_{(A\times B)}},X\right)$ whenever the
$A\times B$ bimodule $M$ is projective on one side. It is clear from the proofs of Section \ref{split} that this
condition can be relaxed, namely it is enough to require the vanishing of the $\Tor$ vector spaces between tensor
powers of the bimodule and the bimodule itself -- we thank Manuel Saorin for stressing this fact. In our situation
the Lemma above shows that $Tor^*_{A\times B}\left(M^{\otimes^p_{(A\times B)}},M\right)=0$ for $p\geq 2$. In order
to show that $Tor^*_{A\times B}\left(M,M\right)=0$, consider a projective resolution of $M$ as a left $B$-module
and extend the action to $A\times B$ letting $A$ act by zero. As in Lemma \ref{tensorzero} tensoring the above
projective resolution by $M$ over $A\times B$ provides a zero complex.

These consideration show that that the cohomology of the columns can be replaced by $\Ext$ vector spaces between
tensor powers of $M$. Much of them vanish using again the Lemma above, finally we obtain the following for the
connecting homomorphism:
\bigskip
$$
\begin{array}{ccc}
&\ \ &\ \ H^{n+1}(A\times B,M)\\ &&\oplus\\ H^n(A\times B,A\times
B)&\stackrel{\delta^{n,0}}{\rightarrow}&\Ext^{n}_{(A\times B)-(A\times B)}(M,M)\\ \oplus &&\oplus\\
\Ext^{n-1}_{(A\times B)-(A\times B)}(M,A\times B)& \stackrel{\delta^{n-1,0}}{\rightarrow}&0
\end{array}
\bigskip
$$ In fact $$\Ext^{n-1}_{(A\times B)-(A\times
B)}(M,A\times B)\ =\ 0\ \ \mbox{and}$$
$$H^{n+1}(A\times B,M)=0.$$ In order to prove this
last assertion, let $e=(1,0)$ and $f=(0,1)$ be the idempotents of the algebra $A\times B$. Note that an $A\times
B$--bimodule $Y$ is the direct sum of four bimodules which can be presented at the vertices of a
square:$$\begin{array}{ll} eYf\ \ &fYf\\ \\ eYe &fYe.
\end{array}$$
For instance $eYf$ is an $A-B$--bimodule and $eYe$ is an $A$--bimodule. We have that {\small$$\begin{array}{lllll}
\Ext_{(A\times B)-(A\times B)}^*(Y, Z)= &\Ext^*_{A-B}(eYf,\ eZf)\ &\oplus\ &\Ext^*_{B-B}(fYf,\ fZf)\ &\oplus\\ \\
&\Ext^*_{A-A}(eYe,\ eZe) &\oplus &\Ext_{B-A}^*(fYe,\ fZe).&
\end{array}$$}
Since the three components $eMf$, $fMf$ and $eMe$ are zero, we obtain $$\Ext^{n-1}_{(A\times B)-(A\times
B)}(M,A\times B)=\Ext^{n-1}_{(A\times B)-(A\times B)}(M,f(A\times B)e),$$ note that $f(A\times B)e=0$. Similarly
we obtain $$H^{n+1}(A\times B,M)=0 \mbox{ since } H^{n+1}(A\times B,M)=\Ext^{n+1}_{(A\times B)-(A\times
B)}(A\times B,M).$$ Moreover the same type of arguments shows that $$H^n(A\times B, A\times B)= H^n(A,A)\oplus
H^n(B,B) \mbox{ and }$$ $$\Ext^n_{(A\times B)-(A\times B)}(M,M)=\Ext^n_{B-A}(M,M).$$ \hfill\qed
\end{Demo}

\begin{obser}
The same result can be derived from the spectral sequence arising from the double complex. Indeed only the first
two columns are non--zero at the first level, and the vector spaces have to be decomposed as we did above.
\end{obser}

\begin{obser}
If $B=k$ and $M$ is any right $A$--module, we obtain Happel's long exact sequence \cite{ha}: $$
\begin{array}{cccccccc}
0&\rightarrow &0&\rightarrow&H^0(T,T)&\rightarrow&H^0(A,A)\oplus k&\rightarrow\\
&&\End_AM&\rightarrow&H^1(T,T)&\rightarrow&H^1(A,A)&\rightarrow\\
&&\Ext^1_A(M,M)&\rightarrow&H^2(T,T)&\rightarrow&H^2(A,A)&\rightarrow\\ &&\Ext^2_A(M,M)&\rightarrow&\cdots
\end{array}$$

\end{obser}

\begin{prop}
The connecting homomorphism of the cohomology long exact sequence for a triangular matrix algebra is given by
$$\delta^nf=1_M\smile f \mbox{ for } f\in H^n(A,A) \mbox{ and }$$
$$\delta^ng=(-1)^{n+1}g\smile 1_M \mbox{ for }g\in
H^n(B,B).$$
\end{prop}
The proof follows from the general description of $\delta^{p,q}$ given in Theorem \ref{smile}.

\footnotesize \noindent C.C.:
\\D\'epartement de Math\'ematiques,
 Universit\'e de Montpellier
2,  \\F--34095 Montpellier cedex 5, France. \\{\tt Claude.Cibils@math.univ-montp2.fr}

\vskip3mm \noindent E.M.:\\ Departamento de Matem\'atica, Universidade de S\~ao Paulo,\\ Caixa Postal 66.281\\
S\~ao Paulo -- SP, 05315--970, Brasil.\\ {\tt enmarcos@ime.usp.br}

\vskip3mm \noindent M.J.R:\\ Departamento de Matem\'atica, Universidad Nacional del Sur,\\Av. Alem 1253\\8000
Bah\'\i a Blanca, Argentina.\\ {\tt mredondo@criba.edu.ar}

\vskip3mm \noindent A.S.:\\ Departamento de Matem\'atica, FCEyN, Universidad de Buenos Aires\\ Pabell\'on I --
Ciudad Universitaria\\ 1428 -- Buenos Aires, Argentina.\\ {\tt asolotar@dm.uba.ar} \vskip3mm \noindent January 21,
2002


\begin{thebibliography}{99}

\bibitem{ibra}
Assem, I.  {\it Alg{\`e}bres et modules}. Enseignement des Math{\'e}matiques.  Les Presses de l'Universit{\'e}
d'Ottawa, Masson, 1997.


\bibitem{asmape} Assem, I., Marcos, E. N.,
 de la Pe\~na, J.A. The Simple connectedness
 of a tame tilted algebra.
J. Algebra 237 (2001) 647--656.

\bibitem{begu} Bendiffalah, B., Guin, D. Cohomologie des morphismes.
 Comm. Algebra 26 (1998) 3939--3951. 18G60


\bibitem{buli} Buchweitz, R.-O., Liu, S. Artin Algebras
with Loops but no Outer Derivations. math.RA/9907008



\bibitem{caei} Cartan, H.; Eilenberg, S.
 {\it Homological algebra}. Princeton University Press, Princeton,
  N. J., 1956

\bibitem{ci} Cibils, C.
Tensor Hochschild homology and cohomology. {\it Interactions between ring theory and representations
  of algebras}
(Murcia), 35--51,
 Lecture Notes in Pure and Appl. Math., 210, Dekker,
 New York, 2000.

\bibitem{ge} Gerstenhaber, M.
On the deformation of rings and algebras. Ann. of Math. 79 (1964), 59--103.

\bibitem{gms} Green, E.L., Marcos, E.N., Snashall, N.
The Hochschild Cohomology Ring of a one point extension, to appear.



\bibitem{ha}
Happel, D. Hochschild cohomology of finite--dimen\-sio\-nal algebras. S\'eminaire d'alg\`ebre Paul Dubreuil et
Marie--Paule Malliavin, Lect. Notes Math. {\bf 1404}, 108--126, 1989.

\bibitem{ha2}
Happel, D. Hochschild cohomology of Auslander algebras. Topics in algebra, Part 1 (Warsaw, 1988), 303--310, Banach
Center Publ., 26, Part 1, PWN, Warsaw, 1990.


\bibitem{mape}  Mart\'{\i}nez-Villa, R., de la Pe\~na,  J. A.
 The universal cover of a quiver with relations, J. Pure
 Appl. Algebra 30 (1983), 277-292.


\bibitem{mipl} Michelena, S., Platzeck, M.I. Hochschild cohomology of triangular matrix
algebras.  J. Algebra 233 (2000), 502--525.


\bibitem{mcc} McCleary, J. {\it User's guide to
 spectral sequences}. Mathematics Lecture Series, 12.
  Publish or Perish, Inc., Wilmington, Del., 1985.


\bibitem{lo}  Loday, J.L. {\it Cyclic homology}.
Grundlehren der Mathematischen Wissenschaften, 301.
 Springer--Verlag, Berlin,
1998.

\bibitem{mj} Redondo, M.J. Universal Galois Coverings of Selfinjective Algebras by
 Repetitive Algebras and Hochschild Cohomology, preprint, 2000.

\bibitem{ma} Saor\'\i n, M. Automorphism groups of trivial extensions,
  J. Pure Appl. Algebra 166 (2002), 285--305.

\bibitem{sk} Skowro\'nski, A. Simply connected
algebras and Hochschild cohomologies, Proc. ICRA IV (Ottawa, 1992), Can. Math. Soc. Conf. Proc. Vol. 14, 431-447,
1993.

\bibitem{we}
Weibel, C.A. {\it An introduction to homological algebra}. Cambridge Studies in Advanced Mathematics, 38.
Cambridge University Press, Cambridge, 1994.

\end{thebibliography}
\end{document}